\magnification=\magstep1
\input amstex
\documentstyle{amsppt}
\hoffset=.25truein
\hsize=6truein
\vsize=8.75truein

\topmatter
\title
A Strong Abhyankar-Moh Theorem and Criterion of Embedded Line
\endtitle
\abstract
    The condition of plane polynomial curve to be a line in well-known Abhyankar-Moh Theorem is replaced by weaker ones. 
A criterion of embedded line is obtained from this strong theorem: Two polynomials can generate the entire polynomial ring iff their derivatives can be generated.
\endabstract
\author
 Yansong Xu
\endauthor
\address
\endaddress
\email
yansong\_xu\@yahoo.com
\endemail
\endtopmatter
\document

\heading    1. Introduction
\endheading

Famous Abhyankar-Moh Theorem [1][2] states that for a field $k$ of characteristic zero, if $f(z)$ and $g(z)$ are polynomials and the polynomial ring $k[f(z),g(z)]=k[z]$, then either $\deg f(z)$  divides $\deg g(z)$  or $\deg g(z)$  divides $\deg f(z)$. But to require the considered polynomial curve to be a line at beginning is too strong and limits the applications of the Theorem. We find that as long as  $\deg f(z)-c$  and $\deg g(z)-c$ are in the degree semigroup of polynomial ring $k[f(z),g(z)]$, where positive integer $c \le \min(\deg f(z), \deg g(z))$, we have either $\deg f(z)$  divides $\deg g(z)$  or $\deg g(z)$  divides $\deg f(z)$. Therefore we  call it Strong Abhyankar-Moh Theorem. Using this strong theorem, we get a criterion for a polynomial plane curve to be an embedded line.

\heading    2. Planar Semigroups
\endheading

Following [5], we first define characteristic $\delta$-sequence and planar semigroup.

\proclaim{Definition 2.1}
Let $\delta=(\delta_0, \delta_1,\cdots,\delta_{h}) \; (h\ge 1)$ be a sequence of $h+1$ natural numbers. And let  $d_k=\gcd(\delta_0, \delta_1,\cdots,\delta_{k-1}) \; (1 \le k\le h+1)$. Then $\delta$ is called a  characteristic         $\delta$-sequence and the semigroup $\Gamma=\Gamma(\delta)$ is called a planar semigroup if the following conditions are satisfied:
\newline
$\operatorname{(1) }$  $d_1\ge d_2 > d_3 > \cdots > d_h > d_{h+1}=1$,
\newline
$\operatorname{(2) }$  $\delta_k \frac{d_k}{d_{k+1}} \in \Gamma(\delta_0, \delta_1,\cdots,\delta_{k-1})  \; (1 \le k \le h)$,
\newline
$\operatorname{(3) }$  $\delta_k < \delta_{k-1} \frac{d_{k-1}}{d_k} \;  (2 \le k \le h)$.
\endproclaim

The following concept of standard expansion is also from [5].

\proclaim{Definition 2.2}
Let $\delta=(\delta_0, \delta_1,\cdots,\delta_{h}) \; (h\ge 1)$ be a  characteristic $\delta$-sequence and let $s$ be an integer. If 
$$ s=a_0\delta_0+ a_1\delta_1+\cdots+a_h\delta_h \; where\; 0 \le a_i < d_i/d_{i+1}  \; (1 \le i \le h) \tag 2.2.1 $$
then we say that $s$ has standard expansion with respect to $\delta$.
\endproclaim

The following properties of characteristic $\delta$-sequence are well-known:
\proclaim{Lemma 2.3 ([5])}
Let $\delta=(\delta_0, \delta_1,\cdots,\delta_{h}) \; (h\ge 1)$ be a  characteristic $\delta$-sequence. Then for any integer $s$,
there is unique standard expansion (2.2.1). Moreover, 
\newline
$\operatorname{(1) }$ $s \in \Gamma(\delta) \text{ iff } a_0 \ge 0$.
\newline
$\operatorname{(2) }$  If $d_i | s$, then $a_j=0$ for $i \le j \le h$.
\endproclaim

We consider standard expansions of two integers and generalize (2) of above Lemma 2.3.

\proclaim{Proposition 2.4}
Let $\delta=(\delta_0, \delta_1,\cdots,\delta_{h}) \; (h\ge 1)$ be a  characteristic $\delta$-sequence.
\newline
$\operatorname{(1) }$  Let $s_1$ and $s_2$ be two integers with standard expansions
$$ s_k=a_{k0}\delta_0+ a_{k1}\delta_1+\cdots+a_{kh}\delta_h\; where \; 0 \le a_{ki} < d_i/d_{i+1} \; (1 \le i \le h) \tag 2.4.1 $$
for $k=1, 2$. Then we have that
 $$ d_i | (s_1 - s_2)  \iff  a_{1j}=a_{2j} \text{ for } i \le j \le h\tag 2.4.2 $$
$\operatorname{(2) }$ If there exists positive integer $c \le \min(\delta_0,\delta_1)$ such that $\delta_0-c \in \Gamma(\delta)$ and $\delta_1-c \in \Gamma(\delta)$, then we have either $\delta_0 | \delta_1$ or $\delta_1 | \delta_0$.

\endproclaim

\demo{Proof} 
(1) $\Longrightarrow $  Since $d_h | d_i$, we have
$$ s_1-s_2=(a_{10}-a_{20})\delta_0+(a_{11}-a_{21})\delta_1+\cdots+(a_{1h}-a_{2h})\delta_h\ \equiv 0 \mod d_h\tag{2.4.3} $$
As we have that
$$ (a_{10}-a_{20})\delta_0+(a_{11}-a_{21})\delta_1+\cdots+(a_{1,h-1}-a_{2,h-1})\delta_{h-1}\ \equiv 0 \mod d_h $$
 $$(d_h/d_{h+1},\delta_h/d_{h+1})=1$$
and
 $$|a_{1h}-a_{2h}|<d_h/d_{h+1}$$
we conclude 
$$a_{1h}=a_{2h}$$
from (2.4.3).
Similarly we can prove that
$$a_{1,h-1}=a_{2,h-1}, \cdots , a_{1i}=a_{2i}.$$
$\Longleftarrow$ Trivial.
\newline
(2) We let $s_1=\delta_0-c$ and  $s_2=\delta_1-c$ in (1). As $s_1-s_2=(\delta_0-c )-  (\delta_1-c)=\delta_0 -  \delta_1$ is divisible by $d_2$,
we have
$$a_{1i}=a_{2i} \; ( 2 \le i \le h)$$
from (1). Without loss of generality, we suppose $\delta_0<\delta_1$. As $c>0$,  it is easy to see that $a_{10}=0$,  $a_{11}=0$ and $a_{21}=0$. Therefore we have
$$\delta_0-\delta_1=s_1-s_2=0-a_{20}\delta_0$$
This proves $\delta_0|\delta_1$
\qed\enddemo

\heading    3. Polynomial Curves
\endheading

We use results of planar semigroups to study polynomial curves.
\proclaim{Theorem 3.1  (Strong Abhyankar-Moh Theorem)}
 Let  $k$ be a field of characteristic zero and $F(f,g)$  be a plane curve which is defined by polynomials  $f(z)$ and  $g(z)$, here $z$ can be an unfaithful parameter.  Let $m$  and $n$  be the degrees of  $f(z)$ and  $g(z)$  respectively. Assume that there is an integer $c>0$ such that $c\leq\min  (m,n)$  and there are polynomials  $u(z)$ and  $v(z)$   in polynomial ring  $k[f(z),g(z)]$ such that  $\deg u(z)=m-c$ and  $\deg v(z)=n-c$, then we have that either  $m$  divides $n$  or $n$  divides $m$.
\endproclaim

\demo{Proof} 
 First we reduce to faithful parameter case. In fact, if $z$  is not a faithful parameter, from [3,  Theorem 3.3.], there exits  $p=p(z)\in k[z]$ and  $\tilde f, \tilde g \in k[z]$ such that  $f(z)=\tilde f (p(z))$ , $g(z)=\tilde g(p(z))$  and $p$  is a faithful parameter. We note that  $u(z) \in k[f(z),g(z)]=k[\tilde f (p(z)),\tilde g (p(z))]$ if and only if there exists $\tilde u(p) \in k[\tilde f(p), \tilde g(p)]$  such that $u(z)=\tilde u(p(z))$ . We also note that $\deg u(z)=m-c$  if and only if  ${\deg}_p \tilde u(p)= \frac{m}{\deg p(z)}-\frac{c}{\deg p(z)}$.  Therefore we only need to handle faithful parameter case.
\newline

From Abhyankar-Moh semigroup theory ([1][2][4][5]), the degree semigroup $\Gamma(F)$ is generated by a $\delta$-sequence $(\delta_0,\cdots, \delta_h)$ with 
$\delta_0=m$ and $\delta_1=n$. Our theorem is a corollary of (2) of Proposition 2.4.
\qed\enddemo

\subheading{Example 3.2}  Let  $f(z)=z^3$ and $g(z)=z^6+z^2$. The plane curve $(f(z),g(z))$  is not an embedded line and we cannot apply Abhyankar-Moh Theorem. Let  $u(z)=z^2$  and  $v(z)=z^5$. As  $u(z),v(z)\in k[f(z),g(z)]$, we can apply Theorem 3.1. Therefore Theorem 3.1 is strictly stronger than Abhyankar-Moh Theorem.

\proclaim{Theorem 3.3 (Criterion of Embedded Line)} Let $k$  be a field of characteristic zero and  $f(z)$  and $g(z)$  be polynomials.  Then $k[f(z),g(z)]=k[z]$  if and only if  $k[f(z),g(z)]\neq k$  and $f' (z)$ and $g' (z)$  are in the polynomial ring $k[f(z),g(z)]$. 
\endproclaim

\demo{Proof}
  $\Longrightarrow$ Trivial.

  $\Longleftarrow$ There are two cases. 

Case 1: $f(z)$ or $g(z)$ is in the field $k$, say $f(z) \in k$  and $g(z) \notin k$. As $g'(z) \in k[f(z),g(z)]=k[g(z)]$, $g(z)$ is linear in $z$, therefore $k[f(z),g(z)]=k[z]$.

Case 2: Both $f(z)$ and $g(z)$ are not in $k$. As $\deg f^{\prime} (z)=\deg f(z)-1$  and $\deg g^{\prime} (z)=\deg g(z)-1$ , we can apply Theorem 3.1 to get the conclusion that either $\deg f(z)$ divides $\deg g(z)$  or   $\deg g(z)$ divides $\deg f(z)$. 
We induct on the sum of  $\deg f(z)$  and  $\deg g(z)$. Let us say $\deg f(z)\geq \deg g(z)$. From Theorem 3.1, we can write $\deg f(z)=l\deg g(z)$,  here $l$ is a positive integer. Let  $a$ and $b$  be the leading coefficients of $f(z)$  and $g(z)$ respectively. We define $f_1(z)=f(z)-a(b^{-1}g(z))^l$  and  $g_1(z)=g(z)$. It is obvious that either $f_1(z)=0$ or   $\deg f_1(z)< \deg f(z)$. It is easy to verify that we still have that $f_1(z)$  is in the polynomial ring $k[f(z),g(z)]=k[f_1(z),g_1(z)]$. We can continue to apply Theorem 3.1. This process must finish in finite steps, say  $n$, and, say $f_n(z) \in k$  and $g_n(z) \in k[f(z),g(z)]=k[f_n(z),g_n(z)]=k[g_n(z)] $. We reduce to case 1, which is proved.
\qed\enddemo

Using the Criterion of Embedded Line, we can give a new equivalence of plane Jacobian Conjecture.

\proclaim {Corollary 3.4}  Plane Jacobian Conjecture is equivalent to the following:

Let $k$  be a field of characteristic zero and  $f(x,y)$  and $g(x,y)$  be polynomials over $k$ with $\deg_{x,y}f(x,y)=\deg_yf(x,y)$.  Assume that $f_xg_y-f_yg_x \in k^*$, then $f_y(x,y)$ and $g_y(x,y)$  are in polynomial ring $k(x)[f(x,y),g(x,y)]$. 
\endproclaim

\demo{Proof } It is well known that Jacobian Conjecture is equivalent to $\deg f|\deg g$  or $\deg g|\deg f$  if  $f_xg_y-f_yg_x \in k^*$  [4]. The corollary is validated by the Criteria of Embedded Line immediately.
\qed\enddemo

\heading    Acknowledgment
\endheading

The author thanks the reviewer for the suggestion to use standard expansion and its properties, which are summarized in Lemma 2.3, to simplify the proof of Theorem 3.1  from the previous version.

\enddocument